\documentclass[12pt]{article}
\usepackage{amsmath}
\usepackage{amssymb}
\usepackage{geometry} 
\title{The Dirichlet series that generates the M\"obius function is
  the inverse of the Riemann zeta function in the right half of the
  critical strip}
\author{Roupam Ghosh}
\date{\today} 

\begin{document}

\maketitle
\begin{center}
\textit{Abstract:}
\\
In this paper I introduce a criterion for the Riemann hypothesis, and then
using that I prove $\sum_{k=1}^\infty \mu(k)/k^s$ converges for $\Re(s) > \frac{1}{2}$. I use a step function $\nu(x) = 2\{x/2\} - \{x\}$ for the Dirichlet eta function ($\{x\}$ is the fractional part of $x$), which was at the core of my investigations, and hence derive the stated result subsequently.
\end{center}
\newpage
In 1859, Bernhard Riemann showed the existence of a deep relationship between two very different mathematical entities., viz. the zeros of an analytic function and prime numbers.
\\\\
The Riemann Hypothesis is usually stated as, the non-trivial zeros of the Riemann zeta function lie on the line $\Re(s) = \frac{1}{2}$. Although, this is the standard formulation, one of the exciting features of this problem is, it can be formulated in many different and unrelated ways. 
\\\\
The approach I take in this paper is influenced by Beurling's 1955 paper: \textit{A closure problem related to the Riemann zeta function} and B\'aez-Duarte's 2001 paper: \textit{New versions of the Nyman-Beurling criterion for the Riemann Hypothesis}, although it takes a new approach. In this paper I would be studying a simple step function $\nu$ relating it to the Dirichlet eta function $\eta$. I will show how the step function $\nu$ convolves with the M\"obius function $\mu(n)$ and gives a constant, which I think is a new result significant at attacking RH.
\\\\
\\\\
\noindent
\textbf{Theorem 1:}
\textit{For the Dirichlet eta function is defined as, for all $\Re(s) > 0$}
$\eta(s) = \sum_{k = 1}^{\infty} \frac{1}{(2k - 1)^s} - \frac{1}{(2k)^s}$
\textit{We get an equivalent expression in the form for all $\Re(s) > 0$
where $ \nu(x) = 2\left\{ x/2 \right\} -  \left\{ x\right\} $, and the expression being given by
\begin{equation}
\eta(s) = s\int_{1}^{\infty} \nu(x) x^{-s-1} dx
\end{equation}
}
\textbf{Proof:}
A simplification of the integral shall prove this case. We have for all $x \in [1,\infty]$ $\nu(x) = 0$ or $1$. It is not hard to see that $\nu(x) = 0$ whenever $x \in [2k,2k+1)$ and $1$ whenever $x \in [2k-1,2k)$ for all positive integers $k$. Hence, we can write the integral as
$$s\int_{1}^{\infty} \nu(x) x^{-s-1} dx = \sum_{k=1}^{\infty} \int_{2k-1}^{2k} s x^{-s-1} dx$$
Giving us the sum
$$\sum_{k = 1}^{\infty} \frac{1}{(2k - 1)^s} - \frac{1}{(2k)^s}$$
which is nothing but $\eta(s)$. Since we already know that this sum converges for $\Re(s) > 0$, we get our result.
\\\\
Alternatively, we can write equation (1) as
\begin{equation}
\frac{\eta(s)}{s} = \int_{0}^{1} \nu\left(\frac{1}{x}\right) x^{s-1} dx
\end{equation}
\textbf{Discussion:} $\nu(x)$ is a simple function that oscillates between 0 and 1 at every integer. As mentioned in the above theorem, $\nu(x) = 0$ whenever $x \in [2k,2k+1)$ and $1$ whenever $x \in [2k-1,2k)$ for all positive integers $k$. Note that $\nu(x/k) = 0$ if $x < k$. The next theorem relates $\nu$ and M\"obius $\mu$ in a very interesting way. I think this might be a new and interesting result, where relating $\mu$ with the simple oscillating step function gives a constant value of $-1$. I think this result is significant for the proof of RH. 
\noindent
\\\\
\\\\
\textbf{Theorem 2:} $\sum_{k=1}^\infty \mu(k)\nu\left(\frac{x}{k}\right) = \left\{
	\begin{array}{ll}
		0  & \mbox{if } x \in [0,1) \\
		1 & \mbox{if } x \in [1,2)	\\
		-1 & \mbox{if } x \in [2,\infty) 
	\end{array}
\right.$
\\\\
\textbf{Proof:}\\
We have for $0 < \theta \leq 1$
$$ \int_{0}^{\theta} \nu\left(\frac{\theta}{x}\right) x^{s-1} dx = \frac{ \theta^s \eta(s)}{s}$$
Since $\nu(\theta/x) = 0$ where $ \theta < x$, or $\theta = 1$ otherwise, hence $$\int_{\theta}^{1} \nu\left(\frac{\theta}{x}\right) x^{s-1} dx = 0 $$
We get the following important equation,
\begin{equation}
\int_{0}^{1} \nu\left(\frac{\theta}{x}\right) x^{s-1} dx = \frac{ \theta^s \eta(s)}{s}
\end{equation}
If we set $f_\mu(x) = \sum_{k=1}^{\infty} \mu(k) \nu(x/k)$ we get
\begin{equation}
\int_{0}^{1} f_\mu\left(\frac{1}{x}\right) x^{s-1} dx = \frac{\eta(s)}{s}\sum_{k=1}^{\infty} \frac{\mu(k)}{k^s}
\end{equation}
\\\\
\textit{(Justification for the exchange of summation and integral: Let $f_n(x) = \sum_{k=1}^{n} \mu(k)\nu(x/k)$. Now since for $\Re(s) > 1$, we have $\sum_{k=1}^{n} \int_{0}^{\infty} |\mu(k)\nu(x/k) x^{-s-1}| dx \leq \sum_{k=1}^{n} \frac{|\mu(k)|}{k^\sigma}\int_{0}^{\infty} \nu(x) x^{-\sigma-1} dx \leq \frac{\zeta(\sigma)\eta(\sigma)}{\zeta(2\sigma)\sigma} < \infty$. By Fubini-Tonelli's theorem we can say,
$\int_{0}^{\infty} \sum_{k=1}^{\infty} \mu(k)\nu(x/k) x^{-s-1} dx = \sum_{k=1}^{\infty} \int_{0}^{\infty} \mu(k)\nu(x/k) x^{-s-1} dx$).}
\\\\
For $\Re(s) > 1$ we know that $\sum_{k=1}^{\infty} \mu(k)/k^s = 1/\zeta(s)$.
Hence, for $\Re(s) > 1$
\begin{equation}
\int_{0}^{1} f_\mu\left(\frac{1}{x}\right) x^{s-1} dx = \frac{1-2^{1-s}}{s}
\end{equation}
Now, since $f_\mu(x) = 1$ whenever $x \in [1,2)$ giving us for all $\Re(s)  > 1$
\begin{equation}
\int_{0}^{\frac{1}{2}} \left( 1 + f_\mu\left(\frac{1}{x}\right) \right) x^{s-1} dx = 0
\end{equation}
Since, $\sum \mu(k) \nu(x/k)$ is always constant in any given $[n,n+1)$. So due to equation (6) we have for all $x \geq 2$ , $f_\mu(x) = -1$. (Properties of Dirichlet series).
\\\\
\textbf{Discussion:} The derivation is not explained above. Notice, the integral in equation (6) above can be expanded as,
$$\frac{-1-f_\mu(2)}{2^s} + \frac{f_\mu(2) - f_\mu(3)}{3^s} + \frac{f_\mu(3) - f_\mu(4)}{4^s} + ... = 0$$
The uniqueness property implies, $-1 - f_\mu(2) = f_\mu(2) - f_\mu(3) = ... = 0$. Since $f_n(x)$ is a step function, we get the result for $x \geq 2$.
\\\\
\\\\
\textbf{Theorem 3:}
\textit{If $f_n = O(n^{\sigma_0+\epsilon})$ where $f_n = |\sup f_n(x)|$,
then $\sum_{k=1}^\infty \mu(k)/k^s$ converges for all $\Re(s) > \sigma_0$.}
\\\\
\textbf{Proof:}
If we set $f_n(x) = \sum_{k=1}^{n} \mu(k) \nu(x/k)$ we get for $\Re(s) > 0$
\begin{equation}
\int_{0}^{\frac{1}{2}} \left(1 + f_n\left(\frac{1}{x}\right)\right) x^{s-1} dx = \frac{\eta(s)}{s}\sum_{k=1}^{n} \frac{\mu(k)}{k^s} - \frac{1 - 2^{1-s}}{s}
\end{equation}
Changing the integral we get
\begin{equation}
\int_{2}^{\infty} \frac{1+f_n\left(x\right)}{x^{s+1}} dx = \frac{\eta(s)}{s}\sum_{k=1}^{n} \frac{\mu(k)}{k^s} - \frac{1-2^{1-s}}{s}
\end{equation}
Since, $1 + f_n(x) = 0$ for $x \in [2,n)$, it follows from above,
\begin{equation}
\left| \frac{\eta(s)}{s}\sum_{k=1}^n \frac{\mu(k)}{k^s} - \frac{1-2^{1-s}}{s}\right| = \int_{2}^{\infty} \frac{\left|1 + f_n(x)\right|}{x^{\sigma+1}} dx = \int_{n}^{\infty} \frac{\left|1 + f_n(x)\right|}{x^{\sigma+1}} dx \leq \frac{\sup|1+f_n(x)|}{\sigma n^{\sigma}}
\end{equation}
If $f_n = O(n^{\sigma_0+\epsilon})$ then LHS converges to $0$ for $\sigma > \sigma_0$ as $n \to \infty$, 
i.e,
\begin{equation}
\lim_{n\to\infty} \left| \frac{\eta(s)}{s}\sum_{k=1}^n \frac{\mu(k)}{k^s} - \frac{1-2^{1-s}}{s}\right| = 0
\end{equation}
which gives, for $\Re(s) > \sigma_0$
\begin{equation}
\frac{\eta(s)}{s}\sum_{k=1}^\infty \frac{\mu(k)}{k^s} = \frac{1-2^{1-s}}{s}
\end{equation}
Now here $\eta(s) \neq 0$ for $1 > \Re(s) > \sigma_0$. \textit{(This is because of, equation (9) with $\eta(s) = 0$ gives LHS $= \left|\frac{1-2^{1-s}}{s}\right|$. Ignoring the line $\Re(s) = 1$, because it is known that M\"obius converges on that line.)}
\\\\
Since in equation (11), $\eta(s) \neq 0$, therefore $\sum_{k=1}^\infty \mu(k)/k^s$ converges for $\Re(s) > \sigma_0$. This is because, we can see from the expression in equation (11), where the M\"obius sum exactly equals $1/\zeta(s)$, which it wouldn't if it did not converge for $\Re(s) > \sigma_0$.
\\\\
\textbf{Theorem 4:}
\textit{$\sum_{k=1}^\infty \mu(k)/k^s$ converges for $\Re(s) > \frac{1}{2}$
}
\\
\textbf{Proof:}\\
For any $\sigma > 0$, $\sigma \neq 1$ the following is justified geometrically.
\begin{equation}
\sum_{k=1}^n \frac{|\mu(k)|}{k^{\sigma}} \leq \int_1^n \frac{dy}{y^{\sigma}} + 1 = \frac{n^{1-\sigma}}{1 - \sigma} - \frac{1}{1 - \sigma} + 1 = O(n^{1-\sigma})
\end{equation}
For $\sigma = 1$ we shall similarly have for any $\epsilon > 0$, $\sum_{k=1}^n |\mu(k)|/k \leq \log(n) + 1 = O(n^\epsilon)$
\\\\
\textbf{I.} Assume $1 \geq \sigma_0 > 1/2$ such that, $\sum_{k=1}^\infty \mu(k)/k^s$ converges \textbf{only} for $\Re(s) > \sigma_0$.
\\\\
\textbf{II.} Now consider for $\sigma > 0$,
\begin{equation}
\begin{split}
\left| \int_1^\infty \frac{f_n(x)}{x^{\sigma+1}} dx \right| \leq \int_1^\infty \frac{|f_n(x)|}{x^{\sigma+1}} dx &= \int_0^\infty \frac{|f_n(x)|}{x^{\sigma+1}} dx\\ 
&\leq \int_0^\infty \frac{\sum_{k=1}^n |\mu(k)|\nu(x/k)}{x^{\sigma+1}} dx\\ 
&= \frac{\eta(\sigma)}{\sigma}\sum_{k=1}^{n} \frac{|\mu(k)|}{k^\sigma}\\ 
&\leq cn^{1-\sigma}
\end{split}
\end{equation}
\\
Assume $\sigma_0 > \sigma' > 1/2$
\begin{equation}
\frac{1}{n^{\sigma'}}\left| \int_1^\infty \frac{f_n(x)}{x^{\sigma+1}} dx \right| \leq \frac{1}{n^{\sigma'}} \int_1^\infty \frac{|f_n(x)|}{x^{\sigma+1}} dx \leq cn^{1 - \sigma' - \sigma}
\end{equation}
Now as $n \to \infty$, equation (12) $\to 0$ for all $\sigma > 1 - \sigma'$. If we consider $$g_n(x) = |f_n(x)|/n^{\sigma'}$$ then the Mellin transform over $g_n(x)$ can be expressed as a Dirichlet series for $\sigma > 1 - \sigma'$, because $g_n(x)$ is a step function (since $f_n(x)$ is a step function) and the transform converges.
$$D_n = \int_1^\infty \frac{g_n(x)}{x^{\sigma+1}} dx = \sum_{k=1}^\infty \frac{a_{n,k}}{k^\sigma}$$

Since the Dirichlet series $D_n \to 0$ for all $\sigma > 1 - \sigma'$ as $n \to \infty$, hence $a_{n,k} \to 0$ as $n \to \infty$ (by the uniqueness property of Dirichlet series, i.e., the series vanishes identically). Giving us $g_n(1) \to 0$, $g_n(1) - g_n(2) \to 0$, $g_n(2) - g_n(3) \to 0$, and so on... as $n \to \infty$. Hence $g_n(x) \to 0$ as $n \to \infty$, and therefore
$$f_n = |\sup f_n(x)| = o(n^{\sigma'})$$
\\\\
\textbf{III.} Now considering theorem 3 and the result in \textbf{II.}, we get $\sum_{k=1}^\infty \mu(k)/k^s$ converges for \textbf{all} $\Re(s) > \sigma'$.
\\\\
\textbf{IV.} But this contradicts our assumption that $\sum_{k=1}^\infty \mu(k)/k^s$ converges \textbf{only} for $\Re(s) > \sigma_0 > \frac{1}{2}$, since $\sigma_0 > \sigma' > 1/2$. Therefore, we must have $\sum_{k=1}^\infty \mu(k)/k^s$ converges for all $\Re(s) > \frac{1}{2}$, thereby validating the Riemann hypothesis.
\\\\

\end{document}